\documentclass[12pt]{article}
\usepackage{epsfig}
\usepackage{latexsym}
\newtheorem{theorem}{Theorem}[section]

\newtheorem{lemma}{Lemma}[section]
\newtheorem{corollary}{Corollary}[section]

\newenvironment{proof}{
\par
\noindent {\bf Proof.}\rm}%
{\mbox{}\hfill\rule{0.5em}{0.809em}\par}
\begin{document}

\title{On the Order of Countable Graphs}

\author{Jaroslav Ne\v{s}et\v{r}il\thanks{Partially
supported by the Project LN00A056 of the Czech Ministry of Education  and by
GAUK 158 grant.}\\ 
Department of Applied Mathematics\\ and\\ Institute for Theoretical
Computer Science (ITI)\\ 
Charles University\\Malostransk\' e n\' am.25, 11800
Praha 1\\ 
Czech Republic\\
{\small Email: nesetril@kam.ms.mff.cuni.cz}\\
\and {Saharon Shelah}\\
Hebrew University, Jerusalem\thanks{Research supported by the United States-Israel Binational
        Science Foundation. Publication 745}\\
Institute of Mathematics\\
         The Hebrew University of Jerusalem\\
         Jerusalem 91904, Israel\\
         and  Department of Mathematics\\
         Rutgers University\\
         New Brunswick, NJ 08854, USA\\
{\small Email: shelah@math.huji.ac.il}}
\maketitle

\footnote{Math. Subj. Class.03C, 05C, 05E, and 06A. Keywords and
phrases: Density, Partially ordered sets, Rigid graphs, Universal graphs. The authors
wish to thank the Mittag-Leffler Institute, for its  support and pleasant
environment during
September, 2000, when this paper was written.}

\vspace{5cm}
                                                   
\begin{abstract}
A set of graphs is said to be {\em independent} if there is no
homomorphism between distinct graphs from the set.
We consider the existence problems related to the independent sets of
countable graphs. While the maximal size of an independent set of countable
graphs is $2^\omega$ the {\em On Line} problem  
of extending an independent set to a larger independent set is much
harder. We prove here that singletons can be  extended (``partnership theorem'').
While
this is the best possible in general, we give structural conditions which
guarantee independent extensions of larger independent sets.

This is related to universal graphs, rigid graphs (where we solve a
problem posed in \cite{NR}) and to the density problem for countable graphs.

\end{abstract}

\section{Introduction and Statement of Results }

Given graphs $G = (V, E),  G' = (V', E')$ a homomorphism  is
any mapping $V \rightarrow V'$ which preserves all the edges of $G$:

\[
\{x, y\} \in E \Longrightarrow \{f(x), f(y)\} \in E' 
\]

This is briefly denoted by $f: G \rightarrow G'$. The existence of a
homomorphism we indicate by $G \rightarrow G'$ and in the context of 
partially ordered sets this will be also denoted by 
$G \leq G'$. $\leq$ is obviously a quasiorder.

$\leq$ is a very rich quasiorder which have been studied in several context,
see \cite{N4} for a survey of this area. For example it has been shown 
(and this also not difficult  to see) 
 that any poset may be represented by $\leq$; see \cite{PT,N2} for
less easy results in this area.
A particular case is an {\em independent set} of graphs which is can be defined
as an independent set (or antichain) in
this  quasiorder.
Here we are interested in a seemingly easy question:

{\em Independence Problem} (shortly IP):

Given a set $\{G_\iota; \iota \in I\}$ of graphs does there exist a graph $G$
such that $\{G_\iota; \iota \in I\}$ together with $G$ form an independent set
of graphs.

This problem has been solved for finite sets of (finite or infinite) graphs in
\cite{N1,N2}. The general case is much harder and it is relative consistent to
assume the negative solution (this is related to Vop\v enka Axiom, see
\cite{ka,N2}).

In this paper we discuss IP for countable graphs.

We prove the following:

\begin{theorem}
\label{one}
For every countable  graph $G$
the following two statements are equivalent:

$i.$ there exists a countable graph
$G'$ such that
$G$ and $G'$ are independent.

$ii.$ $G$ is not bipartite and it does not contain an infinite complete
subgraph.
\end{theorem}

Both conditions given in $ii.$ are clearly necessary. The
non-bipartite comes from the general (cardinality unrestricted)
independence problem as the only finite exception and  the absence of an
infinite clique is due to the cardinality restriction.

This modestly looking result (which we could call {\em Partnership Theorem}:
non-bipartite countable graphs without $K_\omega$ have independent partners)
has a number of consequences and leads to several interesting
problems. First, we want to mention that the above result (and the IP) is
related to universal graphs.

Let $\cal K$ be a class of graphs. We say that a graph $U \in \cal K$ is 
{\em hom-universal} (with respect to $\cal
K$) if $G \leq U$ for every $G \in \cal K$, \cite{NR}.
 
Note that a graph $U$ may be hom-universal with respect to a class $\cal K$
without being universal (in the usual sense: any graph from $\cal K$
is a subgraph of $U$; see \cite{jo,JO,ra,KP,KMP} for an extensive
literature about
universal graphs). For example the triangle $K_3$ is hom-universal for the class 
$\cal K$ of all $3$-colorable graphs and obviously this class does not 
have a finite universal graph. On the other hand clearly any universal 
graph is also hom-universal.

Let $GRA_\omega$ denotes the class of all countable non-bipartite graphs
without an infinite complete subgraph (which is denoted by $K_\omega$).
It is well known  that the class $GRA_\omega$ does not
have a
universal graph. The same proof actually gives that $GRA_\omega$ has no
hom-universal graph. (Here is a simple proof which we sketch for the
completeness: Suppose that $U$ is hom-universal for $GRA_\omega$. Denote by 
$U \oplus x$ the graph which we obtain from $U$ by addition of a new vertex $x$
joined to all the vertices of $U$. Then there exists $f: U \oplus x \rightarrow
U$. Define the vertices $x_0, x_1, \ldots$ by induction: $x_0 = x, x_{i+1} =
f(x_i)$ It is easy to see that all these vertices form a complete graph
in $U$.)

Theorem \ref{one} is a strengthening of the non-hom-universality of
$GRA_\omega$.  In fact Theorem \ref{one} is best possible in the following
sense:

\begin{corollary}
\label{two}
For a positive integer $t$ the following two statements are equivalent:

$i.$ For every finite set $\{G_1, \ldots, G_t\}$  of graphs from $GRA_\omega$
there exists a graph $G\in  GRA_\omega$ such that $G$ and $G_i$
are independent for all $i = 1, \ldots, t$.

$ii.$ t = 1. 
\end{corollary}

\begin{proof}
There are many examples proving $i. \Rightarrow ii.$. For example consider
the complete graph $K_n$ and let $U_n$ be any universal (and thus also hom-universal) 
countable $K_n$-free
graph ($U_n$ exists by \cite{hen}). Then the set $\{K_n, U_n\}$ cannot be extended to a larger independent
set as every graph $G$ either contains $K_n$ or is homomorphic to $U_n$.

An example for $t > 2$ consists from an independent set of finite graphs $G_0,
\ldots , G_{t-1}$ and a countable graph $U$, $U \not\geq G_i, i =
0,\ldots, t-1$ which is universal for all graphs $G$ satisfying $G \not\geq
G_i, i =
0,\ldots, t-1$. Such a graph exists by \cite{CSS,N3}.
(Note also that  an
analogous result does not hold for infinite sets. To see this let $G_i =
C_{2i+3}$ be the set of all cycles of odd length. Then there is no $G$ which
is independent of all graphs $G_i$.) 
\end{proof}

Theorem \ref{one} is in the finite (or cardinality unrestricted) case also
known as (Sparse) Incomparability Lemma \cite{NR,N4}. We can formulate
this as follows:

\begin{theorem}
\label{three}
For any choice of graphs $G, H$, $G$ non-bipartite, satisfying $G < H, \:H
\not\leq G$ there exists a graph $G'$ such that $G' < H, \: H \not\leq G'$ and
such that $G$ and $G'$ are independent.

If $G$ has a finite chromatic number then $G'$ may be chosen finite.
\end{theorem}

(The last part of Theorem \ref{three} may be seen as follows (sketch):
If $\chi(G) = k$ then take $G''$ with $\chi(G'') > k$ and without cycles $\leq
l$ such that $G$ contains an odd cycle of length $\leq l$. Then $G$ and $G''$
are independent. If $\chi(G')$ is large then also graphs $G$ and $G'' \times H$
are independent, see \cite{N1}.)

We do not know whether Theorem \ref{three}  holds
if all the graphs are supposed
to be countable. Partial results are included in Section 5.

Theorem \ref{one} is also related to the notion of a rigid graph:
A graph $G$ is said to be {\em rigid} if its only homomorphism $G \rightarrow
G$ is the identical mapping.
We shall prove the following:

\begin{theorem}
\label{four}
For any countable graph $G$ not containing $K_\omega$ there exists a countable
rigid graph $G'$ containing $G$ as an induced subgraph.
\end{theorem}

The history of this result goes to \cite{CHKN} (the finite case), to \cite{BN}
(the unrestricted cardinality case), and to \cite{NR} (the optimal chromatic
number for  the finite case). Theorem \ref{four} 
solves an open problem proposed in \cite{NR}.

Finally, let us mention that Theorem \ref{one} is related to the concept of
density.

Given a class $\cal K$ of graphs and two graphs $G_1, G_2 \in \cal K$, $G_1 <
G_2$, we say that the pair $(G_1, G_2)$ is a {\em gap} in $\cal K$ if there is
no
$G \in \cal K$ satisfying $G_1 < G < G_2$. The {\em density problem} for class $\cal K$
is the problem to characterize all gaps in $\cal K$. (If there are a ``few''
gaps then we have a tendency to say that class $\cal K$
is {\em dense}; see \cite{NT,N1,N2}.)

Our theorem has the following corollary:

\begin{corollary}
Any pair $(G, K_\omega)$ fails to be a gap in the class of all countable
graphs.
\end{corollary}

\begin{proof}
Let $G < K_\omega$, $G \in GRA_\omega$, be given. According to Theorem
\ref{one} there exists $G' \in GRA_\omega$ such that $G' \not\rightarrow G$.
Then we have $G < G + G' < K_\omega$.
\end{proof}

Note that we used the  easier part of Theorem
\ref{one}. This is being discussed below and some particular positive
examples
of the density of the class $GRA_\omega$ are stated. However the
characterization of all gaps for the class $GRA_\omega$ remains an open
problem. In the class $GRA_\omega$ there are infinitely many gaps. This is in
a sharp contrast with finite graphs where the trivial gap $(K_1, K_2)$ is the
only gap, see \cite{NT}.

The paper is organized as follows:
In Section 2 we give some no-homomorphism conditions which will aid us in
Section 3 in Proof of Theorem \ref{one}. In Section 4 we define high and low
graphs and show  their relationship to the independence problem.
In Section 5 we prove
Theorem \ref{four}. In Section 6 we give structural conditions which allow us
to prove that certain graphs are high and thus generalize Theorems
\ref{one} and \ref {four} to other graphs $H$ than $K_\omega$. We also modify the proof of
Theorem \ref{one} to this setting. This yields a more direct proof
and allows us (at least in principle) to hunt for partners. We  find
classes of graphs where independent extension property holds.  Section 7
contains some remarks and problems.

\section{Necessary conditions for the existence of a homomorphism}

Given two graphs $G_1, G_2$ it is usually not easy to prove that $G_1
\not\rightarrow G_2$. We shall use the following two basic facts:

Suppose $G_1 \rightarrow G_2$. Then 

$i.$ If $G_1$ contains an odd cycle of length $<l$ then also $G_2$ contains an
odd cycle of length $<l$.

$ii.$ $\chi(G_1) \leq \chi(G_2)$ (where $\chi$ denotes the chromatic number).

To this well known list (which cannot be much more expanded even in the finite
case) we add the rank function which we are going to introduce as follows:

Let $G = (V, E)$ be a graph in $GRA_\omega$. By $K_n$ we denote the complete graph on $n =
\{0, 1,\ldots, n-1\}$. Consider the set $h(K_n, G)$
 of all homomorphisms $K_n \rightarrow G$ and denote by $T^G$ the union of all
the sets $h(K_n, G), n = 1, 2, \ldots$. We think of $T^G$ as a (relational)
tree ordered by the relation  $f
\subseteq g$.

It is clear and well known that 

$i.$ $T^G$ is a relational tree;

$ii.$ $T^G$ has no infinite branches.

$iii.$ We can define ordinal ${\rm rk} (T^G) < \omega_1$ the {\em ordinal rank
function} of $T^G$.

(For completeness recall the definition of the ordinal rank function: 
For a tree $T$ without infinite branches  ${\rm rk}(T)$ is defined as
${\rm sup} \{{\rm rk}(T_\iota)+ 1\}$ over all branches of $T$ at the root.)
Put ${\rm rk}(G) = {\rm rk}(T^G)$. We have the following (perhaps folkloristic):

\begin{lemma}
If $G_1 \leq G_2$ then ${\rm rk}(G_1) \leq {\rm rk}(G_2)$.
\end{lemma}

\begin{proof}
Let $f: G_1 \rightarrow G_2$ be a homomorphism. Then for every $n$ we have a
natural mapping $h(f): h(K_n, G_1) \rightarrow h(K_n, G_2)$ defined by 
$h(f)(g) = f \circ g$. $h(f)$ is a level preserving mapping $T^{G_1}
\rightarrow T^{G_2}$ and thus ${\rm rk}(G_1) \leq {\rm rk}(G_2)$.
\end{proof}

For every ordinal $\alpha < \omega_1$ and graph $G$ on $\omega$ consider the
following undirected graph  $K_\omega^{<\alpha>}$:

the vertices of $K_\omega^{<\alpha>}$: all decreasing sequences of ordinal numbers $< \alpha$;

the edges of $K_\omega^{<\alpha>}$:  all pairs $\{\nu, \mu\}$ satisfying $\nu
\subseteq \mu$, by this symbol we mean the  containment of
sequences $\nu$ and $\mu$ (as initial segments).

One can say that $K_\omega^{<\alpha>}$ is a tree of cliques with the total height
$\alpha$.

The following holds for any $\alpha < \omega_1$:

$i.$ $K_\omega^{<\alpha>} \in GRA_\omega$;

$ii.$ ${\rm rk}(K_\omega^{<\alpha>}) = \alpha$.

$iii.$ $K_\omega^{<{\rm rk}(G)+1>} \not\rightarrow G$.

(This gives yet another proof that there is no countable hom-universal
$K_\omega$-free graph.)

Put $G^0 = K_\omega^{<{\rm rk}(G)+1>}$ and let us look at the statement of Theorem \ref{one}.
We have $G^0 \not\leq G$ (by $iii.$) and thus if also $G \not\leq G^0$ then
we are done.
So we can assume the following situation: $G \leq G^0$ and $G^0
\not\leq G$.
Now if $G^1$ is any graph satisfying $G^1 \not\leq G^0$ then necessarily
$G^1 \not\leq G$ (as otherwise $G^1 \leq G \leq G^0$)
and thus by the same token we can assume $G \leq G^1$.
This strategy of the proof will be followed in the next section.

\section{Proof of Theorem \ref{one}}

We proceed by contradiction: Let $G \in GRA_\omega$ be a graph  which
is comparable to every other graph in $GRA_\omega$. By Theorem \ref{three} the chromatic number
of $G$ is infinite.

We shall construct graphs $G^0, G^1, G^2$ such that
$G^0 \not\leq G$ (and thus $G < G^0$), $G^1 \not\leq G^0$ (thus $G < G^1$) and
$G^2 \not\leq G^1$ (and thus  $G < G^2$). Using a construction
similar to the one of $G^2$, we define a family $\{G_\eta\}$ of graphs
which satisfy $G_\eta \not\leq G^1$ and thus $G < G_\eta$. Then the existence of some $\eta$, such that
$G_\eta < G$ will give rise to a contradiction.

The graph $G^0 = K_\omega^{<{\rm rk}(G)+1>}$ was constructed in the previous section.

{\bf Definition of $G^1$:}

The vertices of $G^1$: $\omega \times 2$.
The edges of $G^1$: all pairs of the form $\{(n, i), (m, i)\}$ where $\lfloor
\sqrt n \rfloor = \lfloor
\sqrt m \rfloor$, $i = 0, 1$ and of the form $\{(n, 0), (m, 1)\}$
where $n < m$.

Thus $G^1$  is a ``half graph'' where the vertices are ``blown up'' by
complete graphs of increasing sizes.

{\bf Claim 1:} $G^1 \not\rightarrow G^0$

\begin{proof} ({\em of Claim 1})

Assume contrary: Let $f: G^1 \rightarrow G^0$ be a homomorphism. As $f$
restricted to each of the complete graphs in each of the sets $\omega \times
\{0\}, \omega \times
\{1\}$ is  monotone we can find an infinite set $X \subset \omega$ such that
the mapping $f$ restricted to the set $X \times
\{0\}$ is injective. The set $Y = \{f(x); x \in X\times\{0\}\}$ is an infinite set in 
$V(G^0) = V(K_\omega^{<{\rm rk}(G)+1>})$. The graph $K_\omega^{<{\rm rk}(G)+1>}$
is defined by 
the tree $T, {\rm rk}(T) = {\rm rk}(G)+1$  and
thus by either  K\"{o}nig Lemma (or Ramsey Theorem) the set $Y$ either contains an infinite
chain (i.e. a complete graph in $G^0$) which is impossible, or $Y$ contains an infinite
independent set in $T$ and thus also in $G^0$. 

So $Y$ are the vertices of a star in  $T$ with center $y$. $y$ is
a function $y: K_n \rightarrow G$. Choose $n_\star \in \omega$ such
that the set $X \cap (n_\star \times \{0\})$ has at least $n+1$ elements.

Now the function $f$ restricted to the set $\{n_\star^2, n_\star^2 +1, \ldots
, n_\star^2 + n_\star + 1\} \times \{1\}$ is injective and if $(i,
1)$ is any vertex of this set then $f(i, 1)$ is connected to all vertices
$f(m,0)$ for $m \in X \cap [0, n_\star]$. This implies that $f(m, 0)
\subset y$ for every  $m \in X \cap [0, n_\star]$. But this is a
contradiction. 
\end{proof}

{\bf Construction of $G^2$:}

The vertices of $V(G^2) = A_0 \cup A_1 \cup A_2$ where $A_0 = \{r\}$, and
$A_1$ and $A_2$ are infinite sets (all three mutually disjoint). The set
$A_1$ is disjoint union of finite complete graphs denoted by $K^1_i$
(isomorphic to $K_i$), $i \in \omega$. The set
$A_2$ is disjoint union of finite complete graphs denoted by $K^2_{x,j}$
(isomorphic to $K_j$), $j \in \omega$. The edges of $G^2$ are the edges of all
indicated complete graphs together with all edges of the form $\{r, x\}, x \in
A_1$ and all pairs of the form $\{x, y\}, x \in A_1, y \in \cup_{j\in \omega}
V(K^2_{x, j})$.

So the graph $G^2$ is a tree of depth $2$ with infinite branching with all its
vertices ``blown up'' by complete graphs of increasing sizes.

{\bf Claim 2:}  $G^2 \not\rightarrow G^1$

\begin{proof}

The proof is easy using the main property of the half graph: all the vertices
of one of its ``parts'' (i.e. of the set $\omega \times \{1\}$ ) have finite
degree.

Assume to the contrary that $f: G^2 \rightarrow G^1$ is a homomorphism (for
$G^1$ we preserve all the above notation).
We shall consider two cases according to the value of $f(r)$.

{\rm Case 1}.
$f(r) = (n,1)$ for some $n \in \omega$. 

But then the subgraph of $G^1$ induced
by the neighborhood $N(n,1)$ of the vertex $(n, 1)$ has a finite chromatic
number (as $(n, 1)$ has finite degree in $G^1$) whereas the neighborhood of
$r$ in the graph $G^2$ has the infinite chromatic number (as this neighborhood
is the disjoint union of complete graphs $K^1_i, i \in \omega$).

{\rm Case 2}.
$f(r) = (n,0)$ for some $n \in \omega$.

By a similar argument as in the Case 1 we see that not all vertices $f(x), x
\in A_1$ can be mapped to the vertices of the set $\omega \times \{0\}$ (as by
the connectivity of the subgraph of $G^2$ formed by $A_0 \cup A_1$ this graph
would be mapped to a finite complete graph. 
Thus let $f(x_1) = (m, 1)$ for an $x_1 \in A_1$. But then the
neighborhood $N(m, 1)$ of $(m, 1)$ in the graph $G^2$ has a finite chromatic
number whereas $x_1$  has infinite chromatic number (in $G^1$).
\end{proof}

Thus we see that $G^2 \not\rightarrow G^1$ and consequently $G \rightarrow
G^2$.
The last example which we shall construct will be a family of  graphs $\{G_\eta\}$.
This has to be treated in a more general framewo{\rm rk} and we do it in a separate
subsection.

\subsection{Tree like graphs}

We consider the following generalization of the above construction of $G^2$:

Let $\cal G$ be an infinite set of finite graphs of the form $G_{j,i}, i, j \in \omega$
which satisfies:

$i.$ $\chi(G_{j, i}) \geq i$;

$ii.$ $G_{j,i}$ does not contain odd cycles of length $\leq j$.

$iii.$ All the graphs are vertex disjoint. 
 
Let $T = (V, E)$ be a graph tree (i.e. we consider just
the successor relation) 
defined as follows: 
$V = A_0 \cup
A_1 \cup A_2$ where $A_0 = \{r\}$,
$A_1 = \omega$ and $A_2 = \omega \times \omega$. 
The edges of $T$ are all edges of the form
$\{(r, i)\}, i \in \omega$ and all pairs of the form $\{i, (i, j)\}, i, j\in \omega$.

Let $\eta: V \rightarrow \omega \times \omega$ be any function.

Define the graph
$G_\eta$ as follows:
The set of vertices of $G_\eta$ is the  union of all graphs $G_{\eta(x)}, x
\in V$. The edges of $G_\eta$ are edges of all graphs $G_{\eta(x)}, x
\in \omega$ together with all edges of the form $\{a, b\}$ where $a \in
G_{\eta(x)}$, $b \in G_{\eta(y)}$ and $\{x, y\} \in E$.

Then we have analogously as in the Claim 2:

{\bf Claim 4:}  Let $\eta: V \rightarrow \omega$ be any function
and let $\eta_1, \eta_2: V \rightarrow \omega$ be defined by $\eta(x) = (\eta_1(x),\eta_2(x))$. If $\eta_2$
is unbounded on $A_1$ and on the subsets of $A_2$ of the form $\{i\} \times \omega, i \in \omega$,
then $G_\eta \not\leq G^1$.

Now,  consider the graph $G$ again. As $\chi(G)$ is infinite denote by
$K$ the minimal number of vertices  of a subgraph $G'$ of $G$ with
chromatic number $5$ (by compactness it is $K$ is finite).
Let $\eta: V \rightarrow \omega$ be any
function which  is
unbounded on $\omega$ and each of the sets $\{i\} \times \omega$,  $i \in
\omega$ and moreover which satisfies $\eta_1(i) \geq K$ for every $i \in \omega$.

It is $G_\eta \not\leq G^1$ by Claim 3. Thus $G \leq G_\eta$. In this situation we
prove the following (and this will conclude the proof of Theorem \ref{one}.

{\bf Claim 5:}  $G \not\leq G_\eta$.

\begin{proof}
Assume contrary: let $f: G \rightarrow G_\eta$. Then the vertices of the
subgraph $G'$ are mapped into a set $\cup_{i \in I} G_{\eta(i)}$ where $I$ is
a finite
subset of $V$. Denote by $G''$ the image of $G'$ in $G_\eta$. Due to the tree
structure of $G_\eta$  we have that $\chi(G'') \leq 2 {\rm max}_{i \in I}
\chi(G'' \cap G_{\eta(i)})$.

As  $\eta(i) \geq K$ and thus all graphs $G'' \cap G_{\eta(i)}$ are
bipartite. This implies $\chi(G'' \cap G_{\eta(i)}) \leq 2$ and finally we get
$\chi(G') \leq \chi(G'') \leq 4$, a contradiction.
\end{proof}

\section{Independent Families}

In a certain sense Theorem \ref{one} captures the difficulty of independent
extension property. The pair $K_3, U_3$ (see  Proof following Theorem
\ref{one} in the Section 1) cannot be extended to a large independent set
{\em because} $U_3$ is a rich graph. This can be made precise.
Towards this end we 
first modify the ordinal rank function for graphs below a given graph $H$. We
return to these results in Section 6.

Let $G, H$ be  infinite graphs. Assume that the vertices of $H$ are ordered
in a sequence of type $\omega$.  We can thus assume that $H$ is a graph on
$\omega$. Denote by $H_n$ the subgraph of $H$ induced on the set $\{0, 1,
\ldots , n-1\}$.

Consider the set $h(H_n, G)$
 of all homomorphisms $H_n \rightarrow G$ and denote by $T^G_H$ the union of
all the sets $h(H_n, G), n = 1, 2, \ldots$. We think of $T^G_H$ as a
(relational) tree ordered by the relation  $f
\subseteq g$. $T^G_H$ is called {\em $H$-valued tree of $G$} (with respect of
a given $\omega$-ordering of $H$).

It is clear that 

$i.$ $T^G_H$ is a (relational) tree;

$ii.$ $T^G_H$ has no infinite branches;

 Thus we can define ordinal ${\rm rk}(T^G_H) < \omega_1$ the {\em ordinal rank
function} of $T^G_H$.

Put ${\rm rk}_H(G) =  {\rm rk}(T^G_H)$ ({\em the ordinal $H$-rank} of $G$). We have then
the  following:

\begin{lemma}
Let $G_1, G_2$ be graphs with $H \not\leq G_1$ and $H \not\leq G_1$. Then $G_1
\leq G_2$ implies ${\rm rk}_H(G_1) \leq {\rm rk}_H(G_2)$. 
\end{lemma}

\begin{proof}
Let $f: G_1 \rightarrow G_2$ be a homomorphism. Then  
for every $n$ we have a natural mapping $h(f): h(H_n, G_1) \rightarrow
h(H_n, G_2)$ defined by  $h(f)(g) = f \circ g$. The mapping $h(f)$ is level
preserving mapping  $T_H^{G_1} \rightarrow T_H^{G_2}$ and thus ${\rm rk}_H(G_1) \leq
{\rm rk}_H(G_2)$.   
\end{proof}

For a countable graph $G$ on $\omega$ and every ordinal $\alpha < \omega_1$ define the
following  graph  $G^{<\alpha>}$:

The vertices of $G^{<\alpha>}$ are all decreasing sequences of ordinal numbers $<
\alpha$; the edges of $G^{<\alpha>}$ are  all pairs $\{\nu, \mu\}$ satisfying
$\nu \subseteq \mu$ and $\{\ell(\nu), \ell(\mu)\} \in E(G)$.
(Recall that $\ell(\nu)$ is length of the sequence $\nu$.)

One can say that $G^{<\alpha>}$ is a tree of copies of $G_n$ ($G_n$ is the
graph induced by $G$ on the set $\{0, 1, \ldots, n-1\}$ with the total height
$\alpha$. (This notation also explains the rather cumbersome notation $K_\omega^{<\alpha>}$.)

We have the following:

\begin{lemma}
\label{eight}
$i.$ $G^{<\alpha>} \leq G$;

$ii.$ If $\alpha \leq \beta$ then also $G^{<\alpha>} \leq G^{<\beta>}$;

$iii.$ $G \leq H$ if and only if $G^{<\alpha>} \leq H$ for every $\alpha <
\omega_1$
\end{lemma}

\begin{proof}
This is easy statement. The existing homomorphisms are canonical
level-preserving homomorphisms. Let us mention just $iii.$:

If $f: G \rightarrow H$ then $G^{<\alpha>} \rightarrow H$ by composition of 
$f$ with the map guaranteed by $i.$. Thus assume $G \not\leq H$ and 
$G^{<\alpha>} \leq H$ for any $\alpha < \omega_1$. In this case  
the ordinal $G$-rank of $H$ is defined and ${\rm rk}_G(H) = \alpha < \omega_1$.
As ${\rm rk}_G(G^{<\alpha+1>}) = \alpha+1 > {\rm rk}_G(H)$ we get a contradiction.
\end{proof}

We say that $G$ is $\alpha-${\em low} if $G \leq  G^{<\alpha>}$. A {\em
low graph} is a graph which is low for some $\alpha < \omega_1$, a graph is {\em high} if
it is not low.

We have the following

\begin{theorem}
\label{seven}
Let $G_1, \ldots, G_t$ be an independent set of  countable connected graphs 
including at least one high graph. Then 
there exists a countable  graph $G$ such that $G, G_1, \ldots, G_t$ is an
independent set.
\end{theorem}

\begin{corollary}
Any finite set of high graphs can be extended to a larger independent set.
\end{corollary}

\begin{proof}
Choose the notation such that the graphs $G_1, \ldots, G_{s-1}$ are low while 
graphs $G_s, \ldots, G_t$ are high (the case $s = 0$ corresponds to the set of all high graphs).

Choose $\alpha < \omega_1$ such that for any $m \in \{s, \ldots , t\}$ and 
$n \in \{1, \dots ,s-1\}$ the graph $G_n^{<\alpha>}$ has no homomorphism to
$G_m$. This is possible as by the high-low assumption for every $m, n$ as
above there  is no homomorphism $G_m \rightarrow G_n$ and thus for some
$\alpha(m, n) < \omega_1$ we have $G_n^{<\alpha(m, n)>} \not\longrightarrow
G_m$. (This also covers the case $s = 0$.) Put $\alpha' = {\rm max}\: \alpha(m,
n)$ and $\alpha = {\rm max}\: {\rm rk}_{G_n^{<\alpha(m, n)>}}G_m$.
 
We define

\[
G = \sum _{i = 0}^{t-s} G_{s+i}^{<\alpha>}
\]

and prove that $G$ is the desired graph.
Fix $n \in \{0,\ldots, s-1\}$ and choose $m \in \{s, \ldots , t\}$ arbitrarily.
Then $G_m^{<\alpha>} \not\rightarrow G_n$ and thus $G \not\rightarrow G_n$ as
claimed.

In the opposite direction for every  $m, n \in \{s, \ldots , t\}$ we have 
$G_m  \not\rightarrow G_m^{<\alpha>}$ by $G_m$ high and 
$G_m  \not\rightarrow G_n^{<\alpha>}$ by the choice of $\alpha$ (i.e. as
$\alpha$ is large enough). As $G_m$ is a connected graph $G_m$ maps to $G$
if and only if it maps to one of the component. Thus $G_m \not\rightarrow G$
and we are done.
\end{proof}

{\em Remark}.

Corollary \ref{seven} shows that we have an extension property providing we
``play'' with high graphs. This is in agreement with the ``random building
blocks'' used in the proofs of universality, see \cite{N2}.

\section{Rigid Graphs}

We prove Theorem \ref{four}

Let $G$ be a countable graph not containing $K_\omega$, we can assume that $G$
is  infinite. In fact we can assume without loss of
generality that every edge of $G$ belongs to a triangle and that $G$ is
connected (we simply consider a graph  which contains $G$ as an induced subgraph).

Let $G_1 \in GRA_\omega$ form independent pair with $G$ ($G_1$ exists by
Theorem \ref{one}). We can assume without loss of generality that also every
edge of $G_1$ belongs to a triangle. For that it is enough to attach to every
edge of $G_1$ a pendant triangle; (as every edge of $G$ belongs to a
triangle) these triangles do not influence the non-existence of homomorphisms
between $G$ and $G_1$. $G_1$ can be also assumed to be connected.

Let $G_0$ be a countable rigid graph without triangles. The existence of $G_0$
follows from the existence of a countable infinite rigid relation (take one
way infinite path on $\omega$ together with arc $(0, 3)$) by replacing every
edge by a finite triangle free rigid graph; see e.g. \cite{PT,NR,N4}.

Let $\mu: V(G) \rightarrow V(G_0)$ and $\nu: V(G_0) \rightarrow V(G_1)$ be 
bijections.
Define the graph $G'$ as the disjoint union of graphs $G, G_0, G_1$ together
with the matchings $\{\{x, \mu(x)\}; x \in V(G)\}$ and $\{\{x, \nu(x)\}; x \in
V(G_0)\}$.

We prove that $G'$ is rigid ($G'$ obviously contains $G$ as an induced subgraph).

Let $f: G' \rightarrow G'$ be a homomorphism. As the matching edges and the
edges of $G_0$ do not lie in a triangle we have either $f(V(G)) \subseteq V(G)$
or $f(V(G)) \subseteq V(G_1)$. However the last possibility fails as $G$ and
$G_1$ are independent. Similarly, we have either $f(V(G_1)) \subseteq V(G_1)$
or $f(V(G)) \subseteq V(G)$ and the last possibility again fails.

Thus we have $f(V(G)) \subseteq V(G)$ and $f(V(G_1)) \subseteq V(G_1)$.
As the  vertices  of $G_0$ are the only vertices joined both to $V(G)$ and
$V(G')$ we have also  $f(V(G_0)) \subseteq V(G_0)$. However $G_0$ is rigid and
thus $f(x) = x$ for every $x \in V(G_0)$. Finally as $G$ and $G_0$, $G_0$ and
$G_1$ are joined by a matching we have that $f(x) = x$ for all $x \in V(G')$.

{\em Remark.}

This ``sandwich construction'' may be the easiest proof of a statement of this
type (compare \cite{PT,CHKN,BN,NR}. This proves also the analogous statement
for every infinite $\kappa$ (also for the finite case) providing that we use
the fact that on every set there exists a rigid relation. This has been proved
in \cite{VHP}, and e.g. \cite {N5}
for a recent easy proof.

\section{Gaps below H}

We say that a gap $G < H$ is a {\em gap below $H$}. In
the introduction we derived from Theorem \ref{one} that there are no gaps below
$K_\omega$. It is well known that finite undirected graphs have no non-trivial
gap (except $K_1 < K_2$), see \cite{W,NT}. Also infinite graphs (with
unrestricted cardinalities) have no non-trivial gaps (\cite{N1}). However note
that  classes of graphs with  bounded cardinality (such as $GRA_\omega$)
may have many non-trivial gaps. For example if $H = K_n$ then let $U_n$  be the
hom-universal  $K_n$-free  universal graph. Consider the graph
$G_n = U_n \times K_n$ (the product here is the categorical product defined by
projection - homomorphisms). Then $G_n < K_n$ and it is easy to see that $G_n$
is also $K_n$-free hom-universal graph (universal for graphs below $K_n$). Now if $G < K_n$ then also $G \leq G_n$
and thus  $(G_n, K_n)$ is a gap (below $K_n$). In fact this holds for other
finite graphs, see \cite{N3}. It seems to be difficult to find gaps formed by infinite graphs only. Here we give some explanation of this difficulty. We use the ordinal $H$-rank function for graphs bellow $H$
which was introduced in Section 4.

It is not necessarily true that $H^{<\alpha>} \in GRA_H$. We defined above
$H$ to be an $\alpha$-{\em low graph} if $H^{<\alpha>} \in GRA_H$.
Here are sufficient conditions for low and high:

For a graph $F$ we say that an  infinite subset $X$ of
$V(F)$ is {\em separated by a subset} $C$ if for any two distinct vertices $x, y$ of
$X$ there is no path $x = x_0, x_1, \ldots , x_t = y$ in $F$ such that none of the
vertices $x_1, \ldots , x_{t-1}$  belong to $C$ (thus possibly $x, y \in
C$).

Recall, that graphs $G$ and $G'$ are said to be {\em hom-equivalent} if $G
\leq G' \leq G$. This is denoted by $G \simeq G'$. 

We say that graph $F$ is
{\em $H$-connected} if no infinite subset $X$ of $V(F)$ is separated by a
subset $C$ such that $C \simeq H'$ for a finite subgraph $H'$ of $H$. $H$ is
said to have {\em finite core} if $H$ is equivalent to its finite subgraph. Any
graph with infinite chromatic number has no finite core (and this is far from
being a necessary condition). The following then  holds:

$iv.$ If $H$ is $H$-connected without a finite core then
$H^{<\alpha>} \in GRA_H$.

\begin{proof}
$H$ is infinite. Let $f: H \rightarrow H^{<\alpha>}$ be a homomorphism. As $H$ is not equivalent to any of 
its finite subgraph there exists an infinite set $X \subset V(H)$ such that $f$ restricted to the set $X$ is injective. 
Then the set $f(X)$ is an infinite subset of $H^{<\alpha>}$ and applying the K\"{o}nig's Lemma 
to the tree structure of $H^{<\alpha>}$ we get that either there is an infinite chain 
(which is impossible as $H^{<\alpha>}$ is $H$-free) or there is an infinite star. 
Its vertices form an independent set which is separated by the finite graph corresponding
to the stem of the star.
\end{proof}

We have the following:

\begin{theorem}
\label{six}

Let $H$ be a $H$-connected graph without a finite core.  Then the
following holds:

$i.$ There is no gap below $H$;

$ii.$ $GRA_H$ has no hom-universal graph;

$iii.$ For every $G < H$ there exists $G' < H$ such that $G$ and $G'$ are
independent (``partners under H'').
\end{theorem}

\begin{proof}

$i.$ is easier. Let $G < H$. Then ${\rm rk}_H(G) = \alpha < \omega_1$. It is
$H^{<\alpha+1>} < H$  and thus $H^{<\alpha+1>} \not\rightarrow G$. Put $G^0 = H^{<\alpha+1>}$ and thus we have $G < G +
G^0 < H$ as needed. The same proof gives $ii.$.

However by Lemma \ref{eight} we also know that there exists $\beta > \alpha$
such that $G \not\leq H^{<\beta>}$. This proves $iii.$
\end{proof}

We give another proof of Theorem \ref{six}$iii.$ which is an extension of 
the proof given in Sections 3 and 4. This proof is more direct and gives us more tools for hunting of partners.

\begin{proof}(of Theorem \ref{six},iii.)
 Let $G < H$ be fixed. We proceede in a complete analogy to the above
proof
of Theorem \ref{one} and we outline the main steps and stress only the
differences. Thus let $G$ be a counterexample. Consider $G^0 = H^{<\alpha+1>}$. We have $G^0 \not\leq G$ and
$G^0 < H$ and thus we have $G < G^0$. As $G^0$ has the tree structure we
can find $G^1$ in a similar way such that $G^1 \not\leq G^0$ and $G^1 < H$. 
Given $G^1$ we then define graphs $G^2$ and $G_\mu$ with $G \leq G^2$ and $G
\leq G_\mu$. However we have to continue (as possibly $\chi(H) \leq 4$) and define
also graph $G^4$ with $G \leq G^4$. This will finally lead to a contradiction.

The details of this process are involved and we need several technical
definitions.

An $H$-{\em partite graph} $(G, c)$ is a graph together with a fixed
homomorphism $c: G \rightarrow H$. The sets $c^{-1}(x)$ are {\rm color classes}
of $(G, c)$. Given two $H$-partite graphs $(G, c)$ and $(G', c')$ the $H$-{\em
join} $(G, c) \bowtie (G', c')$ is the disjoint union of $(G, c)$ and $(G', c')$
together with all edges $\{x, x'\}$ where $x \in V(G)$, $x' \in V(G')$ and 
$\{c(x), c(x')\} \in E(H)$. The graph $(G, c) \bowtie (G', c')$ is again
$H$-partite (with the coloring denoted again by $c$).

Recall, that $H_n$ is the graph $H$ restricted to the set $\{0,\ldots, n-1\}$.
Let $H^0_n$ and $H^1_n$ be copies of $H_n$ so that all the graphs $H^0_n$ and
$H^1_n$, $n \in \omega$ are mutually disjoint. Without loss of generality the
vertices of $V(H^i_n)$ belong to $\omega \times \{i\}, i = 0, 1.$
The graphs $H_n, H^0_n, H^1_n$ are considered as $H$-partite graphs with the
inclusion $H$-coloring.

{\bf Definition of $G^1$:}

The vertices of $G^1$ form the set  $\omega \times 2$.
The edges of $G^1$ are  all pairs of the form $\{(x, i), (y, i)\}$ where $\{x, y\}
\in H^i_n$ for some $n \in \omega$ together with all the edges of the graphs 
$H^0_m \bowtie H^1_n$ where $m \leq n$.

$G^1$ is an $H$-partite graph with $c: G^1 \rightarrow H$ defined as the limit
of all the inclusions $H_n \subset H$.
We can still think of $G^1$ as a suitable blowing of a half graph. What is
important is that the key property of half graphs holds for $G^1$: 
all the vertices in the
class $\omega \times \{1\}$ have finite degree.

{\bf Claim 1:} $G^1 \not\rightarrow G^0$

(Recall that  $G^0 = H^{<\alpha+1>}$.)
Assume contrary, let $f: G^1 \rightarrow G^0$.

As $H$ does not have a finite retract we get (by compactness) that for every
$m$ there exists $n$ such that $H_m \not\leq H_n$. It follows that there
exists an infinite set $X \subset \omega$ such that
the mapping $f$ restricted to the set $X \times
\{0\}$ is injective.
 The set $Y = \{f(x); x \in X\}$ is then an infinite subset of
the tree $H^{<\alpha>}, \alpha = {\rm rk}_H(G)$ which defines the graph $G^0$ and
thus by either  K\"{o}nig Lemma (or Ramsey Theorem) the set $Y$ either include an infinite
chain (i.e. a complete graph in $G^0$) which is impossible, or $Y$ include an infinite
independent set in $H^{<\alpha>}$ and thus also in $G^0$. 
So $Y$ are the vertices of an infinite star in  $T_{\alpha, H}$ with center
$y$. $y$ is in fact an injective  homomorphism $y: H_n \rightarrow H$.
Define the set $C$ by $C = f^{-1}(\{0, \ldots, n-1\})$. Then $C$ separates $X$
while $C \leq H_n$. But this is a
contradiction. 
 
{\bf Construction of $G^2$:}

The vertices of $V(G^2) = A_0 \cup A_1 \cup A_2$ where $A_0 = \{r\}$, and
$A_1$ and $A_2$ are infinite sets (all three mutually disjoint). The set
$A_1$ is disjoint union of  graphs $H_i$ denoted by $H^1_i$
(isomorphic to $H_i$), $i \in \omega$. The set
$A_2$ is disjoint union of  graphs denoted by $H^2_{x,j}$
(isomorphic to $H_j$), $x \in A_1, j \in \omega$. The edges of $G^2$ are the
edges of all indicated  graphs $H^1_i$ and $H^2_{x,j}$ together with all edges of the form
$\{r, x\}, x \in A_1$ and all pairs of the form $\{x, y\}, x \in A_1, y \in
\cup_{j\in \omega} V(H^1_{x, j}),\{c(x), c(y)\} \in E(H))$.

So the graph $G^2$ is a tree of depth $2$ with infinite branching with all its
vertices ``blown up'' by  graphs $H_n$ of increasing sizes, the graph
induced by vertices $V(H^1_i) \cup V(H^2_{x,j}$ is isomorphic to $H^1_i
\bowtie H^2_{x,j}$.

{\bf Claim 2:}  $G^2 \not\rightarrow G^1$

\begin{proof}

Assume to the contrary that $f: G^2 \rightarrow G^1$ is a homomorphism (for
$G^1$ we preserve all the above notation).
We shall consider two cases according to the value of $f(r)$.

{\rm Case 1}.
$f(r) = (n,1)$ for some $n \in \omega$. 

(We proceede similarly as in Case 1 of the proof of Theorem \ref{one}.)
But then the subgraph of $G^1$ induced
by the neighborhood $N(n,1)$ of the vertex $(n, 1)$ can be mapped to a
finite subgraph of $H$ (as $(n, 1)$ has finite degree in $G^1$) whereas the
neighborhood of $r$ in the graph $G^2$ cannot be mapped to finite subset of $H$
(as this neighborhood is the disjoint union of  graphs $H^1_i, i \in
\omega$).

{\rm Case 2}.
$f(r) = (n,0)$ for some $n \in \omega$.

This is similar adaptation of the Case 2 of the proof of Theorem \ref{one}.

\end{proof}

Next we shall define graphs $G_\eta$.
We consider the following generalization of the above construction of $G^2$:

Let $\cal G$ be an infinite set of finite graphs of the form $G_{j,i}$
which satisfies:

$i.$ $G_{j, i} \not\rightarrow H_i$;

$ii.$ $G_{j,i}$ do not contain odd cycles of length $\leq j$.

$iii.$ $G_{j, i} \rightarrow H$ (this homomorphism will be denoted again by
$c$). 

$iv.$ All the graphs are vertex disjoint.

By now it is easy to get such examples, see e.g. \cite{NR,N4}.

Let $T = (V, E)$ be a graph tree (i.e. we consider just
the successor relation) 
defined as follows: 
$V = A_0 \cup
A_1 \cup A_2$ where $A_0 = \{r\}$,
$A_1 = \omega$ and $A_2 = \omega \times \omega$. 
The edges of $T$ are all edges of the form
$\{(r, i)\}, i \in \omega$ and all pairs of the form $\{i, (i, j)\}, i, j\in \omega$.

Let $\eta: V \rightarrow \omega \times \omega$ be any function.

Define the graph
$G_\eta$ as follows:
The set of vertices of $G_\eta$ is the  union of all graphs $G_{\eta(x)}, x
\in V$. The edges of $G_\eta$ are edges of all graphs $G_{\eta(x)}, x
\in \omega$ together with all edges of the form $\{a, b\}$ where $a \in
G_{\eta(x)}$, $b \in G_{\eta(y)}$, $\{x, y\} \in E$ and $\{c(A), c(b)\} \in
E(H)$.

We have analogously as in Claim 2:

{\bf Claim 3}  Let $\eta: V \rightarrow \omega$ be any function which is
unbounded on $\omega$ and each of the sets $\{i\} \times \omega$,  $i \in
\omega$. Then $G_\eta \not\leq G^1$.

Now,  consider the graph $G$ again. We have to distinguish two cases:

{\bf Case 1:} $\chi(H) \geq 5$

In this case we proceede completely analogously as in the proof of
Theorem \ref{one} with the only change that we denote by $K$ the minimal
 number of vertices  of  a subgraph
$G'$ of $G$ such that $G' \not\rightarrow H_i$ and $\chi(H_i) \leq 4$
(by compactness it is $K \in \omega$). 
In this case we derive a contradiction as above.
Leaving this at that we have to consider:

{\bf Case 2:} $\chi(H) < 5$

In this case we have to continue and we introduce one more construction
of the graph $G^4$.

Let $T$ be an infinite binary tree. Explicitely, $V(T)$
denotes the set of all binary sequences ordered by the initial segment
containment. For a sequence $\sigma = (\sigma(0),
\sigma(1), \dots , \sigma(p))$ we put $i(\sigma) = \sum_{i=0}^p 2^{\sigma(i)}$
($i(\sigma)$ is a level-preserving enumeration of vertices of $T$) and
$\ell(\sigma) = {\rm max} \:\{i; \sigma(i) \neq 0\}$ ($\ell(\sigma)$ is the level of
$\sigma$ in $T$).

Assume that the graphs $H_n$ satisfy $H_m < H_n$ and $|V(H_m)| < |V(H_n)|$ for all $m < n$. This can be
assumed without loss of generality as we can consider subset of $\omega$ with
this property.  

Let $F_\sigma, \:\sigma \in V(T)$ be a set of disjoint graphs
with the following properties:

$i.$ $F_\sigma \leq H_{i(\sigma)}$.

$ii.$ $F_\sigma > H_{i(\sigma)-1}$, moreover for every homomorphism $f: F_\sigma \rightarrow H$ satisfying $|f(V(F_\sigma))| < |V(F_\sigma)|$ there exist homomorphisms
$g: F_\sigma \rightarrow H_{i(\sigma)}$ and $h: H_{i(\sigma)} \rightarrow H$
such that $f = h \circ g$ (in the other words each $f$ with a small image factorizes through $H_{i(\sigma)}$).

$iii.$ $F_\sigma$ does not contain  odd cycles of length $\leq k_1$ where $k_1$ denotes the shortest length of an odd cycle in $G$.

$iv.$ In each $F_\sigma$ are given two distinct vertices $x_\sigma$ and $y =
y_\sigma$ such that $\{c(x_\sigma), c(y_\sigma)\} \in E(H)$.

(See \cite{NR,N4}; it suffices to put $F_\sigma = H_{i(\sigma)} \times K$
where $K$ is a graph without short odd cycles with sufficiently large
chromatic number.)

Denote by $G^4$ the disjoint union of graphs $F_\sigma$ with
added edges  of the form $\{x, y\}$ where $x = x_\sigma$ and $y = y_\sigma'$ 
and $\{\sigma, \sigma'\}$ form an edge of $T$.

This concludes the definition of $G^4$. For $G^4$ we define $G^3 = G_\eta$ for
the following function $\eta: A_0 \cup A_1 \cup A_2 \rightarrow \omega \time
\omega$ (see the above definition of the graph $G_\eta$ for general $\eta$):

$\eta(r) = 1, \eta(i) = (i, \sum (|V(F_\sigma)|; \ell(\sigma) < i)),
\eta(i, j) = (j, \sum (|V(F_\sigma)|; \ell(\sigma) < j))$.

This only means we consider graphs with rapidly progressing odd girth.

We know that $G^3 \not\rightarrow G^2$ (for any $\eta$ unbounded on the stars
of the
corresponding tree).

Thus assume that $f: G^4 \rightarrow G^3$.
Due to the tree structure of the graph $G^3$ we see that for each $\sigma \in
V(T)$ the image $f(F_\sigma)$ intersects a finite set of graphs
$G_x, x  \in I \subset A_0 \cup
A_1 \cup A_2$ and due to the tree structure of the graph $G^3$ we see easily that there is a homomorphism $f': F_\sigma \rightarrow H_{i(I)}$ where $i(I)$ is the
maximal index appearing among all $i \in I$ and $(j, i) \in I$ and we arrived to a contradiction.

Thus $G^4 \not\leq G^3$ and consequently also $G \leq G^4$.

As $G^4$ contains odd cycles only in copies of graphs $H_\sigma$ and as 
all these cycles have lengths $> k_0$ we conclude that $G \not\leq G^4$.
\end{proof}

\section{Concluding Remarks}

{\bf 1.}
The problem to characterize gaps below $H$ is not as isolated as it perhaps
seems at the first glance. Put $GRA_H = \{G; G < H\}$. We have the following
easy theorem:

\begin{theorem}
\label{five}
For a countable graphs $H$ the following statements are equivalent:

$i.$ There is no gap below $H$;

$ii.$ For every $G \in GRA_H$ there is $G' \in GRA_H$ such
that $G' \not\leq G$;

$iii.$ For every $H' \in GRA_H$ the class $GRA_H'$ has no
hom-universal graph;

\end{theorem}

Motivated by Theorem \ref{four} one is tempted to include here also the
following condition:

$iv.$ For every $G \in  GRA_H$ there exists $G' \in GRA_H$ such that
$G < G'$ and $G'$ is rigid.

However $iv.$ is false as shown by the following example:

Let $H = K_3$ and let $G$ be the disjoint union of all odd cycles of length
$> 3$. Then any rigid graph $G'$, $G' < H$ which contains $G$ as a subgraph is
necessarily a disconnected graph. Let $\{G'_i; i \in \omega\}$ be all the
components of $G'$. Then $\chi(G_i) = 3$ for every $i \in \omega$ and thus let
$G_i$ contains an odd cycle $C_{\ell(i)}$ of length $\ell(i)$. Let $G_j$ be
the component which maps to  $C_{\ell(i)}$ (as a component of $G$). Clearly
$i \neq j$ and thus $G_j \rightarrow G_i$, a contradiction.

Note also that the above Theorem \ref{five} is true for any fixed infinite
cardinality.

{\bf 2.}
We say that a set {\cal G} of countable graphs is {\em maximal} (or {\em
unextendable}) if there is no graph $G \not\in \cal G$ such that $G$ is
independent
to all $G' \in \cal G$.

$\{K_\omega\}$ is maximal but there are other maximal families.
For example $\{K_k\} \cup \{G; G\;  \mbox{finite and}\;\chi(G) > k\}$ is a
maximal set and more generally for every finite graph $H$ the following is a
maximal set:

\[
\{H\} \cup \{G; G\;\mbox{finite and}\; G > H\}
\]

Corollary \ref{two} implies existence of finite maximal sets.

The characterization of maximal sets seems to be difficult problem related to
{\em duality theorems}, see \cite{NT2}. However no maximal set is known which
consists from infinite graphs only.


\begin{thebibliography}{99}

\bibitem{CSS} G. Cherlin, S. Shelah, N. Shi: Universal Graphs with forbidden
subgraphs and
algebraic closure, Advances in Appl. Math. 22 (1999), 454-491.

\bibitem{CS} G. Cherlin, N. Shi: Graphs ommitting a finite set of cycles, J. Graph Th. 21 (1996), 351-355.
 
\bibitem{CHKN}
V. Chv\' atal, P. Hell, L. Ku\v cera, J. Ne\v set\v ril: Every finite graph is
a full subgraph of a rigid graph, J. Comb. Th., 11, 3 (1971), 184-286

\bibitem{BN} L. Babai, J. Ne\v set\v ril: High chromatic rigid graphs II.,
Annals of Discrete Math. 15 (1982), 55-61.

\bibitem{FK} Z. F\"{u}redi, P. Komj\' ath: On the existence of countable
universal graph, J. Graph. Th. 25 (1997), 53-58

\bibitem{FK1} Z. F\"{u}redi, P. Komj\' ath: Nonexistence of universal graphs without some trees, Combinatorica 17 (1997), 163-171.

\bibitem{hen} C. Ward Henson: A family of countable homogeneous graphs, Pacific J. Math. 38 (1971), 69-83.

\bibitem{JO} B. J\' onsson:
Universal Relational Systems, Math. Scand. 4 (1956), 193-208.

\bibitem{jo} J. B. Johnston: Universal infinite partially
orderd sets, Proc. Amer. Math. Soc. 7 (1956), 507-514.

\bibitem{ka}
A. Kanamori: The higher infinite. Large cardinals in set theory from their
beginnings, Perspectives in Math. Logic, Springer Verlag, Berlin 1994

\bibitem{KMP} P. Komj\' ath, A. Mekler, J. Pach: Some Universal Graphs, Israel J. Math. 64 (1998), 158-168.

\bibitem{KP} P. Komj\' ath, J. Pach: Universal graphs without large bipartite
subgraphs, Mathematika 31 (1984), 282-290.


\bibitem{N1} J. Ne\v set\v ril: The Homomorphism Structure of Classes of
Graphs, Combinatorics, Probab. and Comp. 8 (1999), 177-184 

\bibitem{N2}  J. Ne\v set\v ril: The Coloring
Poset and its On - Line Universality, KAM Series 2000-458.
\bibitem{N3}  J. Ne\v set\v ril: Universality and Homomorphism Order (in preparation).
\bibitem{N4}  J. Ne\v set\v ril: Aspects of Structural Combinatorics (Graph
Homomorphisms and their Use), Taiwanese J. Math.3, 4 (1999), 381-424.

\bibitem{N5}  J. Ne\v set\v ril:  A rigid graph
for every set, J. Graph. Th. 39, 2(2002), 108-110. 

\bibitem{NP} J. Ne\v set\v ril, A. Pultr: On Classes of Relations and Graphs
determined by Subobjects and Factorobjects, Discrete Math. 22 (1978), 237-300.

\bibitem{NR} J. Ne\v set\v ril, V. R\" odl: Chromatically Optimal Rigid Graphs,
J. Comb. Th. B, 46 (1989), 133-141.

\bibitem{NT} J. Ne\v set\v ril, C. Tardif: Density. 
In: Contemporary Trends in Discrete Mathematics (R.L.Graham, J.
Kratochv\'{\i}l, J. Ne\v set\v ril, F.S. Roberts, eds.), AMS, 1999, pp. 229-237.

\bibitem{NT2} J. Ne\v set\v ril, C. Tardif: Duality Theorems for Finite
Structures
(characterizing gaps and good characterizatios). J. Comb. Th. B
80 (2000), 80-97.

\bibitem{PT} A. Pultr, V. Trnkov\' a: Combinatorial, Algebraic and
Topological 
Representations of Groups, Semigroups and Categories, North Holland, 1980.

\bibitem{ra} R. Rado: Universal Graphs and Universal Functions, Acta Arith.
9 (1964), 331-340.

\bibitem{W} E. Welzl: Color families are dense, J. Theoret. Comp. Sci. 17 (1982), 29-41.

\bibitem{VHP} P. Vop\v enka, A. Pultr, Z.Hedrl\'{\i}n: A rigid relation
exists on any set,
Comm. Math. Univ. Carol. 6(1965), 149-155

\end{thebibliography}
\end{document}